\documentclass[11pt]{article}

\usepackage[ngerman,french,english]{babel}

\usepackage{pdfsync}

\usepackage{amsmath}                                                        
\numberwithin{equation}{section}
\usepackage{graphicx}                                                       
\usepackage{subfigure}
\usepackage{enumitem}                                                    
\usepackage{hyperref}
\usepackage{verbatim}
\usepackage{xcolor}

\usepackage{amssymb}                                                        
\usepackage[mathscr]{eucal}                                             
\usepackage{cancel}                                                             
\usepackage[normalem]{ulem}                                                                 
\usepackage{pstricks}
\usepackage{rotating}
\usepackage{lscape}
\usepackage[paperwidth=8.5in,paperheight=11in,top=1.00in, bottom=1.00in, left=1.00in, right=1.00in]{geometry}
\usepackage{mathtools}                                                      
\mathtoolsset{showonlyrefs=true}                                    
\usepackage{fixltx2e,amsmath}                                           
\MakeRobust{\eqref}

\usepackage{dsfont}

\usepackage{mathdots}
\usepackage{amsthm}                                                             
\allowdisplaybreaks

\theoremstyle{plain}
\newtheorem{theorem}{Theorem}
\numberwithin{theorem}{section}

\newtheorem{lemma}[theorem]{Lemma}                              
\newtheorem{proposition}[theorem]{Proposition}

\theoremstyle{definition}
\newtheorem{definition}[theorem]{Definition}

\newtheorem{notation}[theorem]{Notation}
\newtheorem{remark}[theorem]{Remark}
\newtheorem{assumption}[theorem]{Assumption}

\def \xy {1}
\def \y {{\eta}}
\def \s {{\sigma}}
\def \g {{\gamma}}
\def \a {{\alpha}}

\def \R {\mathbb{R}}
\def \p {\partial}

\def \EE {\mathsf{E}}
\def \PP {\mathsf{P}}

\def \caratt {{\mathds{1}}}

\newcommand\e{\varepsilon}

\renewcommand\d{\delta}

\renewcommand{\red}[1]{\textcolor{red}{#1}}

\def \Hol {C^{\alpha}}

\def \rn  {{\mathbb {R}}^{N}}
\def \R  {{\mathbb {R}}}
\def \x {{\xi}}
\def \g {{\gamma}}
\def \e {{\varepsilon}}

\def \n {{\nu}}
\def \m {{\mu}}
\def \y {{\eta}}

\def \p {{\partial}}
\def \a {{\alpha}}
\def \O {{\Omega}}

\def \d {{\delta}}

\def \caratt {{\mathds{1}}}

\def \k {{\kappa}}
\def \c {{\chi}}

\def \o {{\omega}}
\def \a {{\alpha}}

\def \d {{\delta}}

\def \G {\Ga}

\def \Ga {{\Gamma}}

\def \It\^o {It\^o }
\def \s {{\sigma}}

\def \w {{\omega}}
\def \R {{\mathbb {R}}}
\def \N {{\mathbb {N}}}

\def \x {{\xi}}
\def \e {{\varepsilon}}

\def \r {{\varrho}}

\def \n {{\nu}}
\def \v {{\nu}}
\def \m {{\mu}}
\def \y {{\eta}}

\def \g {{\gamma}}
\def \O {{\Omega}}
\def \phi {{\varphi}}

\def \tilde {\widetilde}

\def\l {\lambda}

\def \A {\mathcal{A}}

\def \F {\mathcal{F}}

\def \rn {{\mathbb {R}}^{N}}

\def \Ã  {{\`a }}
\def \è {{\`e }}
\def \ò {{\`o }}
\def \ù {{\`u }}

\def \caratt {{\mathds{1}}}

\begin{document}

\title{Existence and uniqueness results for strongly degenerate McKean-Vlasov equations with rough coefficients}

\author{
Andrea Pascucci\thanks{Dipartimento di Matematica, Universit\`a di Bologna, Bologna, Italy.
\textbf{e-mail}: andrea.pascucci@unibo.it} \and Alessio Rondelli\thanks{Dipartimento di
Matematica, Universit\`a di Bologna, Bologna, Italy. \textbf{e-mail}:
alessio.rondelli@studio.unibo.it} \and Alexander Yu Veretennikov\thanks{Institute for Information
Transmission Problems, Moscow, Russian Federation. \textbf{e-mail}: ayv@iitp.ru} }

\date{This version: \today}

\maketitle

\begin{abstract}
We present existence results for weak solutions to a broad class of degenerate McKean-Vlasov
equations with rough coefficients, expanding upon and refining the techniques recently introduced
by the third author. Under certain structural conditions, we also establish results concerning
both weak and strong well-posedness.
\end{abstract}



%
%

\section{Introduction}\label{intro}
We consider the McKean-Vlasov (MKV) equation in $\R^{N}$
\begin{equation}\label{e1}
  dX_{t}=B(t,X_{t},\m_{X_{t}})dt+\Sigma(t,X_{t},\m_{X_{t}})dW_{t}
\end{equation}
where $\m_{X_{t}}$ denotes the law of $X_{t}$ and $W_{t}$ is a $d$-dimensional Brownian motion on
a probability space $(\O,\F,\PP)$. We assume $0\le d\le N$, meaning the equation can be
degenerate, even {\it completely}. The coefficients may have linear growth and minimal regularity
properties: as a key example, our framework encompasses the following MKV-Langevin equation with
{\it measurable} coefficients
\begin{equation}\label{eL}
  \begin{cases}
    dX_{0,t}=X_{1,t}dt, \\
    dX_{1,t}=B_{1}(t,X_{t},\m_{X_{t}})dt
    +\Sigma_{1}(t,X_{t},\m_{X_{t}})dW_{t}.
  \end{cases}
\end{equation}
The classical Langevin model \cite{Langevin} is a specific instance of \eqref{eL}: the solution is
a $2d$-dimensional process $X_{t}=(X_{0,t},X_{1,t})$ describing the dynamics in the phase space of
a system of $d$ particles with position $X_{0,t}$ and velocity $X_{1,t}$ at time $t$. The interest
in measurable coefficients is primarily driven by applications in control problems. In finance,
SDEs of the form \eqref{eL} describe path-dependent contingent claims, such as Asian options or
some local stochastic volatility model (see, for instance, \cite{MR2791231}).

Stochastic differential equations dependent on distributions have been extensively researched
since McKean's seminal paper \cite{MR0221595}, building upon Kac's foundation of kinetic theory
\cite{MR0084985}. Genuinely degenerate MKV equations have recently sparked considerable interest:
they are studied in \cite{MR2834721} as alternative approaches to Navier-Stokes equations
turbulent flows, using the results in \cite{DiFrancescoPascucci2}. The martingale problem for
nonlinear kinetic distribution-dependent SDEs with singular drifts is studied in \cite{RusPagl},
\cite{Anh} and \cite{zhang2022second}, which also contain numerous additional references. Recent
results for stable driven MKV equations are proved in \cite{Chaudru1}, \cite{Chaudru2} and
\cite{Hao}. Applications in the calibration of local stochastic volatility models in finance and
nonlinear filtering problems are discussed in \cite{MR3224295}, \cite{MR4152640} and
\cite{MR4563699}, while recent work on degenerate quantile-dependent SDEs can be found in
\cite{MR3072241} and \cite{MR4713517}. Rather general existence results for the martingale problem
for \eqref{e1} were proved in \cite{MR762085} assuming the continuity of the coefficients and the
existence of Lyapunov functions.

\medskip In order to state our main results, we introduce the assumptions and notations we will adopt throughout the paper. The
subsequent assumptions can be significantly loosened, and expansions to coefficients of wider
generality are outlined in Remark \ref{r2}.
\begin{assumption}[\bf Structural assumptions]\label{H1}
The coefficients of \eqref{e1} are of the form
\begin{equation}\label{e2}
 B(t,x,\m)=\int_{\R^{N}}b(t,x,y)\m(dy),\qquad \Sigma(t,x,\m)=\int_{\R^{N}}\s(t,x,y)\m(dy),\qquad
 t\in[0,T],\ x\in\R^{N},
\end{equation}
where $b$ is a $\R^{N}$-valued Borel measurable function and $\s$ is a $\R^{N\times d}$-valued
Borel measurable function of the form
\begin{equation}\label{e13}
 \s=\begin{pmatrix}
     0 \\
     \s_{1} \\
   \end{pmatrix}.
\end{equation}
In \eqref{e13}, $0$ is a null $(N-d)\times d$-block, while the block $\s_{1}$ is uniformly
positive definite, that is
\begin{equation}\label{e20}
  \langle\s_{1}\x,\x\rangle\ge\l|\x|^{2},\qquad \x\in\R^{d},
\end{equation}
for some positive constant $\l$.
\end{assumption}
Later it will be convenient to rewrite \eqref{e1} in a more explicit form, separating the
degenerate from the diffusive part of the SDE as in \eqref{eL}.
\begin{notation}\label{n2} We set $x=(x_{0},x_{1})\in\R^{N-d}\times\R^{d}$ and call $x_{0}$ and $x_{1}$ the {\it
degenerate} and {\it non-degenerate} components of $x\in\R^{N}$, respectively.
\end{notation}
Setting $X_{t}=\left(X_{0,t},X_{1,t}\right)$, equation \eqref{e1} reads
\begin{equation}\label{e1split}
  \begin{cases}
    dX_{0,t}=B_{0}(t,X_{t},\m_{X_{t}})dt,\\
    dX_{1,t}=B_{1}(t,X_{t},\m_{X_{t}})dt+\Sigma_{1}(t,X_{t},\m_{X_{t}})dW_{t}.
  \end{cases}
\end{equation}
As per Notation \ref{n2}, we say that $B_{0}$ is the {\it degenerate} drift coefficient and
$B_{1},\Sigma_{1}$ are the {\it non-degenerate} drift and diffusion coefficients respectively,
which are defined according to \eqref{e2}:
\begin{align}
  B_{i}(t,x,\m_{X_{t}})&=\int_{\R^{N}}b_{i}(t,x,y)\m_{X_{t}}(dy)=\EE\left[b_{i}(t,x,X_{t})\right],\qquad i=0,1,\\ \label{e12}
  \Sigma_{1}(t,x,\m_{X_{t}})&=\int_{\R^{N}}\s_{1}(t,x,y)\m_{X_{t}}(dy)=\EE\left[\s_{1}(t,x,X_{t})\right],\qquad
 t\in[0,T],\ x\in\R^{N},
\end{align}
with $\s_{1}$ as in \eqref{e13}.
\begin{assumption}\label{H3}
The coefficients $b=(b_{0},b_{1})$ and $\s_{1}$ satisfy the linear growth condition
\begin{equation}\label{e23}
  |b(t,x,y)|+|\s_{1}(t,x,y)|\le \mathbf{c}(1+|x|+|y|),\qquad t\in[0,T],\, x,y\in\R^{N},
\end{equation}
for some positive constant $\mathbf{c}$. Moreover, $b=b(t,x_{0},x_{1},y_{0},y_{1})$ and
$\s_{1}=\s_{1}(t,x_{0},x_{1},y_{0},y_{1})$ are continuous functions of the degenerate variables
$x_{0},y_{0}\in\R^{N-d}$,  uniformly in $t\in[0,T]$ and $x_{1},y_{1}\in\R^{d}$.
\end{assumption}
Assumption \ref{H3} imposes the continuity of the coefficients solely with respect to the {\it
degenerate} variables $x_{0},y_{0}$. This encompasses both ends of the spectrum: from
deterministic equations with {\it continuous} coefficients (as per the classic Cauchy-Peano
existence theorem, but see also \cite{MR2487734}), to non-degenerate stochastic equations with
{\it measurable} coefficients (as per the well-known results by Krylov \cite{MR270462} and
\cite{MR2723141}, Chap.2 Sec.6).

{\begin{definition}[\bf Weak solution]\label{d1} A weak solution on $[0,T]$ to the MKV equation
with coefficients $b$ and $\s_{1}$, is a pair $(X_{t},W_{t})_{t\in[0,T]}$  of stochastic
processes, defined on a filtered probability space $(\O,\F,\PP,(\F_{t})_{t\in[0,T]})$, such that:
\begin{itemize}
  \item[i)] $W$ is a $d$-dimensional Brownian motion;
  \item[ii)] $X$ is a continuous and adapted process such that
\begin{equation}\label{ae11}
  \int_{0}^{T}\EE\left[|X_{t}|\right]^{2}dt<\infty;
\end{equation}
  \item[iii)] almost surely we have
\begin{equation}\label{ae12}
  X_{t}=X_{0}+\int_{0}^{t}\left(\int_{\R^{N}}b(s,X_{s},y)\m_{X_{s}}(dy)\right)ds
  +\int_{0}^{t}\left(\int_{\R^{N}}\s_{1}(s,X_{s},y)\m_{X_{s}}(dy)\right)dW_{s},\qquad t\in[0,T].
\end{equation}
\end{itemize}
\end{definition}
\begin{remark}\label{ra1}
Under Assumption \ref{H3}, condition \eqref{ae11} ensures that the stochastic integral in \eqref{ae12} is well-defined. In
fact, we have
\begin{equation}
 \int_{0}^{T}\left|\int_{\R^{N}}\s_{1}(t,X_{t},y)\m_{X_{t}}(dy)\right|^{2}dt\le
 \mathbf{c}^{2}\int_{0}^{T}\left(1+|X_{t}|+\EE\left[|X_{t}|\right]\right)^{2}dt
\end{equation}
which is finite almost surely, by \eqref{ae11} and since $X$ is continuous. Notice that for non-MKV equations condition \eqref{ae11} is
redundant since it follows from standard a priori $L^{p}$ estimates for solutions of SDEs (cf.,
for instance, \cite{PascucciSC}, Theor. 14.5.2).
\end{remark}}
Our first result is the following
\begin{theorem}[\bf Existence]\label{t1}
Let $\m_{0}$ be a distribution on $\R^{N}$ with finite fourth moment and $T>0$. Under Assumptions
\ref{H1} and
\ref{H3}, 
a weak solution $(X,W)$ on $[0,T]$ to the MKV equation with coefficients $b,\s_{1}$ and with
initial distribution $\m_{X_{0}}=\m_{0}$, exists and verifies
\begin{align}\label{e4aa}
  \EE\bigg[\sup_{0\le t\le T}|X_{t}|^{4}\bigg]\le
  C(1+\EE\left[|X_{0}|^{4}\right]),\qquad 
  \sup_{0\le t,s\le T\atop |t-s|\le h}\EE\left[|X_{t}-X_{s}|^{4}\right]\le C
  h^{2},
\end{align}
for some positive constant $C$.
\end{theorem}


Existence results were recently established by the third author of this paper in
\cite{Veretennikov} under stronger conditions, specifically requiring the coefficients to be
bounded (excluding the prototype MKV-Langevin equation \eqref{eL}), the coefficient $b_{0}$ to be
continuous with respect to all variables except $t$, and the diffusion matrix to be symmetric. The
main tools used in the proof are Krylov's bounds \cite{MR2723141}, Skorokhod's technique of weak
convergence \cite{MR185620}, and Nisio's approach to SDEs in \cite{MR334328} and \cite{MR4100228}.
In addition to extending the original arguments of \cite{Veretennikov}, in Proposition \ref{pp1}
we clarify certain technical aspects of the construction of the so-called $\e$-net, a crucial tool
in the proof of Theorem \ref{t1}: see, in particular, Remark \ref{rrr1}.

\begin{remark}\label{r2}
The techniques employed in the proof of Theorem \ref{t1} are quite likely to be applicable even in
the case of coefficients of the form
\begin{equation}\label{e2b}
 \bar{B}(t,x,\m)=\Phi(t,x,B(t,x,\m)),\qquad \bar{\Sigma}_{1}(t,x,\m)=\Psi(t,x,\Sigma_{1}(t,x,\m)),\qquad
 t\in[0,T],\ x\in\R^{N},
\end{equation}
where $B,\Sigma_{1}$ in \eqref{e2}-\eqref{e12} satisfy Assumptions \ref{H1}-\ref{H3}, and
  $$\Phi:[0,T]\times\R^{N}\times\R^{N}\longrightarrow \R^{N},\qquad
  \Psi:[0,T]\times\R^{N}\times\R^{d\times d}\longrightarrow \R^{d\times d},$$
are measurable functions, locally Lipschitz continuous in the second and third argument, uniformly
in $t\in[0,T]$, and $\bar{\Sigma}_{1}$ satisfies the non-degeneracy condition \eqref{e20bis}. Also
the assumption regarding the finiteness of the fourth moment of the initial distribution $\m_{0}$
can be relaxed: the proof remains valid without modification as long as $\m_{0}$ has a finite
$(2+\e)$-moment, for some $\e>0$.
\end{remark}

\medskip
In the second part of the paper, we establish the well-posedness of the MKV equation
\eqref{e1split} under the specific
condition that the coefficients $B_{0}$ and $\Sigma_{1}$ are {independent of the law}. We 
consider the MKV equation
\begin{equation}\label{e1bb}
  \begin{cases}
    dX_{0,t}=B_{0}(t,X_{t})dt,\\
    dX_{1,t}=B_{1}(t,X_{t},\m_{X_{t}})dt+\Sigma_{1}(t,X_{t})dW_{t},
  \end{cases}
\end{equation}
with initial datum $\y$. We also consider the standard (i.e. non-MKV) equation
{\begin{equation}\label{e1nonMKV}
  \begin{cases}
    d\hat{X}_{0,t}=B_{0}(t,\hat{X}_{t})dt,\\
    d\hat{X}_{1,t}=B_{1}(t,\hat{X}_{t},\m_{t})dt+\Sigma_{1}(t,\hat{X}_{t})dW_{t},
  \end{cases}
\end{equation}
}where $\m=(\m_{t})_{t\in[0,T]}$ is the flow of marginals of a given continuous, adapted process.
Equation \eqref{e1nonMKV} is a linearized version of \eqref{e1bb}.

{\begin{theorem}[\bf Uniqueness]\label{at1}
Suppose the structural Assumption \ref{H1} holds and the coefficients are Borel measurable functions
of the form
  $$b=(b_{0}(t,x),b_{1}(t,x,y)),\qquad \s_{1}=\s_{1}(t,x),\qquad t\in[0,T],\, x,y\in\R^{N}.$$
Assume also that:
\begin{itemize}
  \item[i)] $b$ satisfies the linear growth condition
  \begin{equation}\label{e23c}
  |b(t,x,y)|\le \mathbf{c}(1+|x|),\qquad t\in[0,T],\, x,y\in\R^{N};
\end{equation}
  \item[ii)] the matrix $\s_{1}$ is bounded
and uniformly positive definite, that is
\begin{equation}\label{e20b}
  \l^{-1}|\x|^{2}\le\langle\s_{1}(t,x)\x,\x\rangle\le\l|\x|^{2},\qquad t\in[0,T],\,
  x\in\R^{N},\ \x\in\R^{d},
\end{equation}
for some positive constant $\l$;
  \item[iii)] the initial datum satisfies
\begin{equation}\label{e23cd}
  \EE\left[e^{r |\y|^{2}}\right]<\infty
\end{equation}
for some positive constant $r$.
\end{itemize}
If, for any fixed flow $\m$ the solution of \eqref{e1nonMKV}, with initial datum $\y$, is weakly
(resp. strongly) unique, then also the solution of \eqref{e1bb} is weakly (resp. strongly) unique.
\end{theorem}}
\begin{remark}\label{rr1}
In Theorem \ref{at1}, when $b_{1}$ is bounded, a less restrictive condition than \eqref{e23cd} is
sufficient, such as $\EE\left[|\y|^{2}\right]<\infty$.
\end{remark}
Notice that Theorem \ref{at1} does not require any regularity conditions for the coefficients, which may be merely measurable.
On the other hand, one of the most stringent assumptions 
is the weak (or strong) well-posedness of
the non-MKV equation \eqref{e1nonMKV}. 
To illustrate the potential range of applications of the theorem, we provide a far from complete
list of known results on well-posedness for \eqref{e1nonMKV}. We distinguish between the {\it non-degenerate} case
(i.e., when $d = N$) and the {\it degenerate} case (i.e., when $d < N$):
\begin{itemize}
\item \cite{MR185620}, \cite{MR2190038}, \cite{MR270462,MR0341607} prove weak well-posedness
 for non-degenerate SDE with bounded and measurable coefficients
  (see also Sect. 2.6 in \cite{MR2723141}  and Theorems 6.1.7 and 7.2.1 in \cite{MR2190038}.
  In \cite{MR0336813}, \cite{MR568986} strong well-posedness is established if the drift is bounded and H\"older
  continuous, and the diffusion coefficient is bounded and Lipschitz continuous;
\item several results for degenerate Langevin-type SDEs are available in the literature: for SDEs with H\"older
continuous coefficients, weak well-posedness is studied in \cite{MR3758337} and is a straightforward consequence of the
results in \cite{MR4660246} and \cite{DiFrancescoPascucci2} where an optimal notion of intrinsic
H\"older continuity is considered. Well-posedness of the martingale problem for SDEs with
measurable diffusion coefficients is established in \cite{zhang2022second}. Strong well-posedness
is achieved when the coefficients are bounded and H\"older continuous, provided the H\"older
exponents meet certain thresholds, as proved in \cite{MR4035024} and \cite{MR4498506}. For
instance, in the case
  $B_{0}(t,X_{0,t},X_{1,t})=X_{1,t}$ the drift $B_{1}$ must be at least $2/3$-H\"older
continuous w.r.t. the degenerate component; it is currently not known if this threshold is indeed
optimal. The author of \cite{MR688921} establishes sufficient conditions for the strong
well-posedness of SDEs with nonsmooth drift under conditions of possible degeneracy of the
diffusion with respect to some of the variables.
\end{itemize}
Due to the technique employed, which relies on Girsanov's theorem, in this paper
we prove well-posedness only for MKV equations whose diffusion coefficient is
independent of the law. In
future research, we plan to revisit the study of well-posedness under more general assumptions.
One approach we are considering involves combining techniques from \cite{MR3072241} and
\cite{MR4035024} for non-degenerate equations with the parametrix method used in
\cite{DiFrancescoPascucci2} and \cite{MR4660246} for the linear case. This combined approach could
allow us to exploit the non-Euclidean intrinsic geometric structures induced by degenerate kinetic
equations, as detailed in \cite{MR3429628}.

\medskip The paper is organized into two sections, which respectively contain the proofs of the
existence and uniqueness theorems. Throughout the paper, $C$ denotes a positive constant that may
vary from line to line.

\section{Existence: proof of Theorem \ref{t1}}
\subsection{Step 1: regularization}
We approximate the coefficients through a sequence of functions $b^{n},\s^{n}$ that are bounded,
Lipschitz continuous and satisfy Assumptions \ref{H1} and \ref{H3}. Specifically, they satisfy the
non-degeneracy condition \eqref{e20} and the linear growth condition \eqref{e23} with the
constants $\l,\mathbf{c}$ {\it independent of $n$}. This can be achieved by taking, for $f=b,\s$,
the convolution $f^{n}:=f_{n}\ast\psi_{n}$ where $\psi_{n}=\psi_{n}(t,x,y)$ is a standard
mollifier supported in the ball of radius $\frac{1}{n}$, and
 $$f_{n}(t,x_{0},x_{1},y_{0},y_{1}):=f(t,\c_{n}(x_{0}),\c_{n}(x_{1}),\c_{n}(y_{0}),\c_{n}(y_{1}))$$
where, with a slight abuse of notation, we set
    $$\c_{n}(\x)=
    \begin{cases}
    \x & \text{if } |\x|\le n, \\
    n\frac{\x}{|\x|}& \text{otherwise},
   \end{cases}$$
both for $\x\in\R^{d}$ and $\x\in\R^{N-d}$. Clearly, $f_{n}\equiv f$ on $[0,T]\times Q_{n}$ where
\begin{equation}\label{e25}
  Q_{n}:=\{(x_{0},x_{1},y_{0},y_{1})\in\R^{2N}\mid \max\{|x_{0}|,|x_{1}|,|y_{0}|,|y_{1}|\}\le  n\}.
\end{equation}
{To ensure the convolution is well-defined, we extend the domains of $b$ and $\sigma$ to include
$t<0$. We define $b(t,\cdot)=0$ and $\sigma(t,\cdot)=I_n$ for $t<0$. The extension of $\sigma$ in
this manner is designed to maintain its positive definiteness.} We denote by $B^{n},\Sigma^{n}$
the integrals corresponding to $b^{n},\s^{n}$, as defined in \eqref{e2}.

The equation with regularized coefficients is strongly well-posed. More precisely, let $W^{x}$ be
a $d$-dimensional Brownian motion defined on a filtered probability space $(\O,\F,\PP,\F_{t})$. A
fixed point argument (see, for instance, \cite{MR1108185} or \cite{MR3752669}, Theor. 4.21) shows
that, for any $n\in\N$ and $\y_{x}\in L^{p}(\O,m\F_{0},\PP)$ with $p\ge 2$, the MKV equation
\begin{equation}\label{e1b}
  dX^{n}_{t}=B^{n}(t,X^{n}_{t},\m_{X^{n}_{t}})dt+\Sigma^{n}(t,X^{n}_{t},\m_{X^{n}_{t}})dW^{x}_{t},\qquad
  X^{n}_{0}=\y_{x},
\end{equation}
admits  a unique solution in the space $\mathbb{S}^{p}$ of continuous processes, adapted to the
augmented filtration $\F^{\y_{x},W^{x}}$ generated by $\y_{x}$ and $W^{x}$, and such that
\begin{equation}\label{ae7}
  \EE\left[\sup_{0\le t\le T}|X^{n}_{t}|^{p}\right]<\infty.
\end{equation}
Next, as a consequence of the linear growth condition \eqref{e23}, we derive some estimates that
are uniform in $n$. Although the proof is quite standard, the authors were unable to locate an
appropriate reference; therefore we provide a sketch of the proof for completeness.
\begin{proposition}\label{ap1} Let $\y_{x}\in
L^{p}(\O,m\F_{0},\PP)$ with $p\ge 2$. For any $n\in\N$ the solution $X^{n}\in \mathbb{S}^{p}$ of
\eqref{e1b} satisfies the estimates
\begin{align}\label{e4}
  \EE\left[\sup_{0\le t\le T}|X^{n}_{t}|^{p}\right]&\le
  C(1+\EE\left[|\y_{x}|^{p}\right]),\\ \label{e5}
  \sup_{0\le t,s\le T\atop |t-s|\le h}\EE\left[|X^{n}_{t}-X^{n}_{s}|^{p}\right]&\le C
  h^{p/2}.
\end{align}
where the constant $C$ depends on $T,N$ and $p$, but not on $n$.
\end{proposition}
\proof We prove that, for any $n\in\N$,
\begin{equation}\label{ae6}
  v_{n}(t_1):=\EE\left[\sup_{0\leq t\leq t_1}|X^{n}_t|^p\right]\leq C\left(1 + \EE\left[|\y_{x}|^p\right]\right)\left(1 +
 \int_0^{t_1}v_{n}(s)d s\right),\qquad 0\le t_{1}\le T,
\end{equation}
where $C$ is a constant that depends only on $p,N,T$ and $\mathbf{c}$ in \eqref{e23}: in
particular, $C$ is independent of $n$. Since $v_{n}(T)$ is finite by \eqref{ae7}, we can directly
apply Gr\"onwall's lemma to deduce \eqref{e4} from \eqref{ae6}.

In order to prove \eqref{ae6}, we notice that by \eqref{e23} we have
\begin{equation}\label{ae5}
  |B^{n}(t,X^{n}_{t},\m_{X^{n}_{t}})|+|\Sigma^{n}(t,X^{n}_{t},\m_{X^{n}_{t}})| 
  \le \mathbf{c}
  \int_{\R^{N}}(1+|X^{n}_{t}|+|x|)\m_{X^{n}_{t}}(dx)=
  \mathbf{c}(1+|X^{n}_{t}|+\EE\left[|X^{n}_{t}|\right]).
\end{equation}
Applying H\"older and Burkholder-Davis-Gundy inequalities to \eqref{e1b}, for some constant
$c_{0}$ dependent only on $p$ and $N$, we get
\begin{align}
  v_n(t_1) &\leq 3^{p-1}\left(\EE\left[|\y_x|^p\right] +
  t_1^{p-1}\int_0^{t_1}\EE\left[|B^{n}(s,X^{n}_s,\mu_{X^{n}_s})|^p\right]d s
  +c_{0} t_1^{\frac{p-2}{2}}\int_0^{t_1}\EE\left[|\Sigma^{n}(s,X^{n}_s,\mu_{X^{n}_s})|^p\right] d s\right)\le
\intertext{(by \eqref{ae5}, for some constant $c_{1}$ dependent only on $p,N,T$ and $\mathbf{c}$)}
  &\leq c_{1} 
  \left(\EE\left[|\y_x|^p\right] +
  \int_0^{t_{1}}
  \EE\left[1+|X^{n}_s|^p+\EE\left[|X^{n}_s|^p\right]\right]ds\right)
\end{align}
which yields \eqref{ae6}.
The proof of \eqref{e5} is based on similar arguments and is omitted.
\endproof


\subsection{Step 2: linearization
}
By possibly enlarging the probability space $(\O,\F,\PP)$, 
we may assume the existence of a second Brownian motion $W^{y}$, independent of $W^{x}$, and a
random variable $\y_{y}$ equal in law to $\y_{x}$ and independent of $\y_{x}$. We consider the
sequence $\left(Y^{n}\right)_{n\in\N}$ of strong solutions of the SDE \eqref{e1b} where $W^{y}$
substitutes $W^{x}$ and $\y_{y}$ replaces $\y_{x}$. By uniqueness, $Y^{n}$ is equal in law to
$X^{n}$ and, being adapted to $\F^{\y_{y},W^{y}}$, is independent of $W^{x}$ (hence, also
independent of $X^{m}$ for any $m\in\N$).

Notice that $X^{n}$ solves the standard (i.e. non-MKV) SDE
\begin{equation}\label{e1c}
  dX^{n}_{t}=B^{n}(t,X^{n}_{t},\m_{Y^{n}_{t}})dt+\Sigma^{n}(t,X^{n}_{t},\m_{Y^{n}_{t}})dW^{x}_{t},\qquad
  X^{n}_{0}=\y_{x}.
\end{equation}
Clearly, in equation \eqref{e1c} we can swap the roles of $X^{n}$ and $Y^{n}$: thus, we have
\begin{equation}\label{e1cy}
  dY^{n}_{t}=B^{n}(t,Y^{n}_{t},\m_{X^{n}_{t}})dt+\Sigma^{n}(t,Y^{n}_{t},\m_{X^{n}_{t}})dW^{y}_{t},\qquad
  Y^{n}_{0}=\y_{y},
\end{equation}
and $Y^{n}$ satisfies estimates \eqref{e4}-\eqref{e5} as well.


\subsection{Step 3: limiting process}
For convenience, we introduce the compact notations $Z=(X,Y)$ and $W=(W^{x},W^{y})$.
By virtue of estimates \eqref{e4}-\eqref{e5} and Skorokhod's lemma, there exists a probability
space supporting a sequence of stochastic processes $(\tilde{Z}^{n},\tilde{W}^{n})$ that are equal
in law to $(Z^{n},W)$ and such that, up to a subsequence, converge in probability
\begin{equation}\label{e7}
  (\tilde{Z}_{t}^{n},\tilde{W}^{n}_{t})\xrightarrow[n\to\infty]{\PP}
  (\tilde{Z}_{t}^{\infty},\tilde{W}^{\infty}_{t})
\end{equation}
for all $t\in[0,T]$. 

\begin{remark}\label{r1}
The fact that \eqref{e7} holds for any $t\in[0,T]$ implies the weak convergence of
finite-dimensional distributions of the processes $(\tilde{Z}^{n},\tilde{W}^{n})$.
\end{remark}
\begin{lemma}\label{p1}
For any $n\in\N\cup\{\infty\}$, the process $\tilde{Z}^{n}$ satisfies estimates
\eqref{e4}-\eqref{e5} with the constant $C$ independent of $n$:
\begin{align}\label{e4Z}
  \EE\left[\sup_{0\le t\le T}|\tilde{Z}^{n}_{t}|^{4}\right]&\le
  C(1+\EE\left[|\y_{x}|^{4}\right]),\\ \label{e5Z}
  \sup_{0\le t,s\le T\atop |t-s|\le h}\EE\left[|\tilde{Z}^{n}_{t}-\tilde{Z}^{n}_{s}|^{4}\right]&\le C
  h^{2}.
\end{align}
Moreover, $\tilde{Z}^{n}$ admits a continuous modification and $\tilde{W}^{n}$ is a
$2d$-dimensional Brownian motion with respect to the filtration $(\tilde{\F}^{n}_{t})_{t\in[0,T]}$
generated by $(\tilde{Z}^{n},\tilde{W}^{n})$.
\end{lemma}
\proof For $n\in\N$ the thesis follows from \eqref{e4}-\eqref{e5} and the equivalence in law of
$Z^{n}$ and $\tilde{Z}^{n}$. If $n=\infty$, by Fatou's lemma\footnote{Up to sub-sequences, we have
a.s. convergence of $\tilde{X}^{n}_{t}$ to $\tilde{X}^{\infty}_{t}$ for any $t$.} we have
\begin{equation}\label{e5b}
  \sup_{0\le t,s\le T\atop |t-s|\le h}\EE\left[|\tilde{Z}^{\infty}_{t}-\tilde{Z}^{\infty}_{s}|^{4}\right]\le
  \sup_{0\le t,s\le T\atop |t-s|\le h}\liminf_{n\to \infty}
  \EE\left[|\tilde{Z}^{n}_{t}-\tilde{Z}^{n}_{s}|^{4}\right]\le C h^{2}
\end{equation}
which proves \eqref{e5Z}. By Kolmogorov's continuity theorem, $\tilde{Z}^{n}$ admits a continuous
modification for any $n\in\N\cup\{\infty\}$.

Next, we prove \eqref{e4Z} with $n=\infty$. 
Let $\mathcal{D}_{m}$ denote the $m$-th dyadic set of $[0,T]$. Since $(\tilde{Z}^{n}_{t})_{t\in
\mathcal{D}_{m}}$ converges in probability as $n\to\infty$ to $(\tilde{Z}^{\infty}_{t})_{t\in
\mathcal{D}_{m}}$, we also have a.s. convergence up to a subsequence of $\max\limits_{t\in
\mathcal{D}_{m}}|\tilde{Z}^{n}_{t}|^{4}$ to $\max\limits_{t\in
\mathcal{D}_{m}}|\tilde{Z}^{\infty}_{t}|^{4}$: thus, by Fatou lemma and \eqref{e4} we have
\begin{equation}\label{e4bb}
  \EE\left[\max_{t\in\mathcal{D}_{m}}|\tilde{Z}^{\infty}_{t}|^{4}\right]\le
  \liminf_{n\to \infty}\EE\left[\max_{t\in\mathcal{D}_{m}}|\tilde{Z}^{n}_{t}|^{4}\right]\le
  C(1+\EE\left[|\y_{x}|^{4}\right]).
\end{equation}
By Beppo-Levi we can pass in the limit in \eqref{e4bb} as $m\to \infty$. Hence \eqref{e4Z} follows
from the continuity of $\tilde{Z}^{\infty}$.

Now, let $n\in\N$: since $(\tilde{Z}^{n},\tilde{W}^{n})$ is equal in law to $(Z^{n},W^{n})$, then
$\tilde{W}^{n}_{t}-\tilde{W}^{n}_{s}$ is independent of
$(\tilde{Z}^{n}_{r},\tilde{W}^{n}_{r})_{r\in[0,s]}$ (and therefore, of $\tilde{\F}^{n}_{s}$).
Using also that $\tilde{W}^{n}$ is a Gaussian process (because it is equal in law to $W$), this
suffices to prove that $\tilde{W}^{n}$ is a Brownian motion with respect to the filtration
$(\tilde{\F}^{n}_{t})_{t\in[0,T]}$. {This result extends to $n=\infty$ by virtue of Remark
\ref{r1}.}
\endproof

\subsection{Step 4: SDE for $\tilde{X}^{n}_{t}$ with $n\in\N$}
For any $n\in\N\cup\{\infty\}$, $\tilde{Z}^{n}$ is continuous and adapted to
$(\tilde{\F}^{n}_{t})_{t\in[0,T]}$, and therefore progressively measurable: by \eqref{e12} and
Fubini's theorem, also the processes $\Sigma^{n}(t,\tilde{X}^{n}_{t},\m_{\tilde{Y}^{n}_{t}})$ are
progressively measurable and therefore the stochastic integrals
\begin{equation}\label{e8}
  \int_{0}^{t}\Sigma^{n}(s,\tilde{X}^{n}_{s},\m_{\tilde{Y}^{n}_{s}})d\tilde{W}^{x,n}_{s},\qquad
  n\in\N\cup\{\infty\},
\end{equation}
are well-defined.
\begin{lemma}\label{p2}
For any $n\in\N$, $(\tilde{X}^{n},\tilde{W}^{x,n})$ solves the SDE
\begin{equation}\label{e1ct}
  d\tilde{X}^{n}_{t}=B^{n}(t,\tilde{X}^{n}_{t},\m_{\tilde{Y}^{n}_{t}})dt+\Sigma^{n}(t,\tilde{X}^{n}_{t},\m_{\tilde{Y}^{n}_{t}})d\tilde{W}^{x,n}_{t}.
\end{equation}
An analogous result is valid for $(\tilde{Y}^{n},\tilde{W}^{y,n})$.
\end{lemma}
The proof of Lemma \ref{p2} derives from a standard approximation argument (cf. \cite{MR2723141},
p.89) and is postponed to Section \ref{sp1}.

\subsection{Step 5: SDE for $\tilde{X}_{t}^{\infty}$}
A key step in the proof is the following proposition which states that we can pass in the limit in
\eqref{e1ct}. Here we rely on the compactness argument used in \cite{Veretennikov} that allows to
``freeze'' the degenerate component $\tilde{X}^{n}_{0,t}$ in order to apply {\it Krylov's
estimates} to the non-degenerate component $\tilde{X}^{n}_{1,t}$. Since the coefficients have
linear growth at infinity, we employ the broader form of Krylov's bounds provided by Lemma 3 of
\cite{MR4421344}.
\begin{proposition}\label{pp1}
The process $\tilde{X}^{\infty}_{t}$ in \eqref{e7} solves the SDE
  \begin{equation}\label{e15}
  \tilde{X}^{\infty}_{t}=\tilde{X}^{\infty}_{0}+\int_{0}^{t}B(s,\tilde{X}^{\infty}_{s},
  \m_{\tilde{Y}^{\infty}_{s}})ds+\int_{0}^{t}\Sigma(s,\tilde{X}^{\infty}_{s},
  \m_{\tilde{Y}^{\infty}_{s}})d\tilde{W}^{x,\infty}_{s}.
\end{equation}
\end{proposition}
\proof We split the proof in three steps. We prove the $L^{1}$-convergence of the drift and the
convergence in probability of the diffusive part as we take the limit in \eqref{e1ct} as $n$ tends
to infinity.

\medskip\noindent{\bf Step I: the $\e$-net.}\\
Let $\e$ and $\a$ be two fixed positive constants. By Ascoli-Arzela's theorem, for any $h\in\N$ the space 
\begin{equation}\label{rrre1}
  \Hol_{h}:=\{\phi\in C([0,T],\R^{N-d})\mid |\phi(t)|\le h,\, |\phi(t)-\phi(s)|\le h|t-s|^\alpha,\,t,s\in[0,T]\}
\end{equation}
is totally bounded. Hence, there exists a finite collection of functions
$\phi_{1}^{(h)},\dots,\phi^{(h)}_{\k_{\e}}\in\Hol_{h}$, which we call an {\it $\e$-net}, such that
\begin{equation}\label{e22}
 \Hol_{h}=\bigcup_{j=1}^{\k_{\e}}Q_{\e,h}(\phi_{j}^{(h)}),\qquad
  Q_{\e,h}(\phi_{j}^{(h)}):=\{\phi\in\Hol_{h}\mid \sup_{t\in[0,T]}|\phi(t)-\phi^{(h)}_{j}(t)|< \e\}.
\end{equation}
Here $\k_{\e}\in\N$ depends only on $\e$, $\a$, $T$, the dimension $N-d$ and $h\in\N$. Now, assume
$\a<1/2$: by \eqref{e5Z} and Kolmogorov's continuity theorem, $\tilde{Z}^{n}(\o)\in
C^{\a}([0,T],\R^{2N})$ for any $\w\in\O$ and $n\in\N$. Moreover, by \eqref{e4Z}-\eqref{e5Z} and
Markov inequality, there exists $h\in\N$ such that
\begin{equation}\label{e22b}
\PP\left((\tilde{X}_{0,\cdot}^n\in C^{\a}_{h})\cap(\tilde{Y}_{0,\cdot}^n\in C^{\a}_{h})\right)\ge
  1-\frac{\e}{2},\qquad  n\in\N\cup\{\infty\},
\end{equation}
with $h$ dependent on $\e$ but not on $n\in\N\cup\{\infty\}$. Thus, we also have
\begin{equation}\label{e16pre}
 \PP\left(\mathcal{B}_{n,\e}\right)\geq 1-\frac{\e}{2},\qquad
 \mathcal{B}_{n,\e}:=\bigcup_{i,j=1}^{\k_{\e}}\left(\tilde{X}_{0,\cdot}^n\in
 Q_{\e,h}(\phi_{i}^{(h)})\right)\cap\left(\tilde{Y}_{0,\cdot}^n\in Q_{\e,h}(\phi_{j}^{(h)})\right)
\end{equation}
for any $n\in\N\cup\{\infty\}$. Again, notice that $\k_{\e}$ depends on $\e$, $\a$, $T$ and the
dimension $N-d$, but not on $n$.

By estimate \eqref{e4Z}, which is uniform in $n$, for any $\e>0$ there exists $M_{\e}>0$ such that
\begin{equation}\label{e16}
 P(\mathcal{D}_{n,{\e}})\ge 1- \frac{\e}{2},\qquad
 \mathcal{D}_{n,{\e}}:=\bigg(\sup_{0\le t\le T}|\tilde{Z}^{n}_{t}|<M_{\e}\bigg)\cap\bigg(\sup_{0\le t\le
 T}|\tilde{Z}^{\infty}_{t}|<M_{\e}\bigg)
\end{equation}
for any $n\in\N\cup\{\infty\}$.
In conclusion, setting
\begin{equation}\label{e14Q}
  \mathcal{Q}_{n,\e}:=\mathcal{B}_{n,\e}\cap\mathcal{D}_{n,{\e}},
\end{equation}
by \eqref{e16pre} and \eqref{e16} we have
\begin{equation}\label{e16bisp}
  \PP(\mathcal{Q}_{n,\e})\geq 1-\e,\qquad
  \mathcal{Q}_{n,\e}=\bigcup_{i,j=1}^{\k_{\e}}\mathcal{Q}^{i,j}_{n,\e},\qquad
  n\in\N\cup\{\infty\},
\end{equation}
where
\begin{equation}\label{e18}
  \mathcal{Q}^{i,j}_{n,\e}:=\left(\tilde{X}^{n}_{0,\cdot}\in Q_{\e,h}(\phi_{i}^{(h)})\right)
  \cap\left(\tilde{Y}^{n}_{0,\cdot}\in Q_{\e,h}(\phi_{j}^{(h)})\right)\cap
  \mathcal{D}_{n,{\e}},\qquad 1\le i,j\le \k_{\e},
\end{equation}
with $Q_{\e,h}(\phi_{i}^{(h)})$ as in \eqref{e22}. Notice that, without loss of generality, we can
modify the sets $\mathcal{Q}^{i,j}_{n,\e}$ to be mutually disjoint.
\begin{remark}\label{rrr1}
By the uniform estimates \eqref{e4Z}-\eqref{e5Z} and the total boundedness of $C^{\a}_{h}$ in
\eqref{rrre1}, the number $\k_{\e}$ of elements of the $\e$-net is independent of $n$. We
explicitly indicate the dependence of $\k_{\e}$ on $\e$ due to the delicate estimates such as
\eqref{e28b} where $\k_{\e}$ appears as the upper limit in the summation (see \cite{Veretennikov},
page 13 for comparison).
\end{remark}


\medskip\noindent{\bf Step II: $L^{1}$-convergence of the drift.}\\
We prove that, for any arbitrary $\d>0$,
  $$\lim_{n\to \infty}\EE\left[\int_{0}^{t}
  \left|B^{n}(s,\tilde{X}^{n}_{s},\m_{\tilde{Y}^{n}_{s}})-
  B(s,\tilde{X}^{\infty}_{s},\m_{\tilde{Y}^{\infty}_{s}})\right|ds\right]\le\d.
  $$
Let us fix $\d>0$ and notice that, for some $m\in\mathbb{N}$ which depends on $\delta$ and will be
chosen appropriately later, we have
\begin{align}\label{e33}
  \int_{0}^{t}\EE\left[
  \left|B^{n}(s,\tilde{X}^{n}_{s},\m_{\tilde{Y}^{n}_{s}})-
  B(s,\tilde{X}^{\infty}_{s},\m_{\tilde{Y}^{\infty}_{s}})\right|\right]ds\le \mathbf{E}_{1}^{n,m}+\mathbf{E}_{2}^{n,m}+\mathbf{E}_{3}^{m}
\end{align}
where
\begin{align}\label{e36}
  \mathbf{E}_{1}^{n,m}&:=\int_{0}^{t}\EE\left[
  \left|B^{n}(s,\tilde{X}^{n}_{s},\m_{\tilde{Y}^{n}_{s}})-
  B^{m}(s,\tilde{X}^{n}_{s},\m_{\tilde{Y}^{n}_{s}})\right|\right]ds,\\
  \mathbf{E}_{2}^{n,m}&:=\int_{0}^{t}\EE\left[
  \left|B^{m}(s,\tilde{X}^{n}_{s},\m_{\tilde{Y}^{n}_{s}})-
  B^{m}(s,\tilde{X}^{\infty}_{s},\m_{\tilde{Y}^{\infty}_{s}})\right|\right]ds,\\ \label{e37}
  \mathbf{E}_{3}^{m}&:=\int_{0}^{t}\EE\left[
  \left|B^{m}(s,\tilde{X}^{\infty}_{s},\m_{\tilde{Y}^{\infty}_{s}})-
  B(s,\tilde{X}^{\infty}_{s},\m_{\tilde{Y}^{\infty}_{s}})\right|\right]ds.
\end{align}
\begin{itemize}
\item[$\diamond$] {\bf The term $\mathbf{E}_{1}^{n,m}$.}\\
We consider $\mathcal{Q}_{n,\e}$ in \eqref{e14Q}: the value of $\e$ will be selected later
depending on $\d$. We have
\begin{align}
  \mathbf{E}_{1}^{n,m}&\le\int_{0}^{t}\EE\left[
  \EE\left[\left|b^{n}(s,x,\tilde{Y}^{n}_{s})-
  b^{m}(s,x,\tilde{Y}^{n}_{s})\right|
  \right]\!\!\raisebox{-0.11ex}{$\Big\vert$}_{x=\tilde{X}^{n}_{s}}\right]ds=\label{e32}
\intertext{(by the freezing lemma, being $\tilde{X}^{n}_{s}$ independent of $\tilde{Y}^{n}_{s}$)}
  &=\int_{0}^{t}\EE\left[\left|b^{n}(s,\tilde{X}^{n}_{s},\tilde{Y}^{n}_{s})-
  b^{m}(s,\tilde{X}^{n}_{s},\tilde{Y}^{n}_{s})\right|\right]ds\le\\
&\le\int_{0}^{t}\EE\left[\left(\mathds{1}_{\mathcal{Q}_{n,\e}}+\mathds{1}_{\Omega\setminus\mathcal{Q}_{n,\e}}\right)\left|b^{n}(s,\tilde{X}^{n}_{s},\tilde{Y}^{n}_{s})-
  b^{m}(s,\tilde{X}^{n}_{s},\tilde{Y}^{n}_{s})\right|\right]ds\le\\
\intertext{(by the linear growth of $b$ uniformly in $n$, the fact that
$\PP\left(\Omega\setminus\mathcal{Q}_{n,\e}\right)\le\e$ and \eqref{e4Z})} \label{e28}
  &\le \e\,CT(1+\EE\left[|\y_{x}|^{4}\right])+I^{n,m}
\end{align}
where, by \eqref{e18},
\begin{align}
 I^{n,m}&:=
  \sum_{1\le i,j\le \k_{\e}}\int_{0}^{t}\EE\left[\caratt_{\mathcal{Q}^{i,j}_{n,\e}}\left|b^{n}(s,\tilde{X}^{n}_{s},\tilde{Y}^{n}_{s})-
  b^{m}(s,\tilde{X}^{n}_{s},\tilde{Y}^{n}_{s})\right|\right]ds\\
  &\le\sum_{1\le i,j\le \k_{\e}}(I^{n,m,i,j}_{1}+I^{n,m,i,j}_{2}+I^{n,m,i,j}_{3})
\end{align}
with $\mathcal{Q}^{i,j}_{n,\e}$ as in \eqref{e18} and
\begin{align}
  I^{n,m,i,j}_{1}&=\int_{0}^{t}\EE\left[\caratt_{\mathcal{Q}^{i,j}_{n,\e}}
  \left|b^{n}(s,\tilde{X}^{n}_{0,s},\tilde{X}^{n}_{1,s},\tilde{Y}^{n}_{0,s},\tilde{Y}^{n}_{1,s})-
  b^{n}(s,\phi_{i}(s),\tilde{X}^{n}_{1,s},\phi_{j}(s),\tilde{Y}^{n}_{1,s})\right|\right]ds,\\
  I^{n,m,i,j}_{2}&=\int_{0}^{t}\EE\left[\caratt_{\mathcal{Q}^{i,j}_{n,\e}}\left|b^{n}(s,\phi_{i}(s),\tilde{X}^{n}_{1,s},\phi_{j}(s),\tilde{Y}^{n}_{1,s})-
  b^{m}(s,\phi_{i}(s),\tilde{X}^{n}_{1,s},\phi_{j}(s),\tilde{Y}^{n}_{1,s})\right|\right]ds,\\
  I^{n,m,i,j}_{3}&=\int_{0}^{t}\EE\left[\caratt_{\mathcal{Q}^{i,j}_{n,\e}}\left|
  b^{m}(s,\phi_{i}(s),\tilde{X}^{n}_{1,s},\phi_{j}(s),\tilde{Y}^{n}_{1,s})-
  b^{m}(s,\tilde{X}^{n}_{0,s},\tilde{X}^{n}_{1,s},\tilde{Y}^{n}_{0,s},\tilde{Y}^{n}_{1,s})\right|\right]ds.
\end{align}
By Assumption \ref{H3}, $b^{n}(s,\cdot,\cdot)$ and $b^{m}(s,\cdot,\cdot)$ are uniformly continuous
on $Q_{M_{\e}}$ in \eqref{e25}, with respect to the degenerate variables $x_{0},y_{0}$, uniformly
in $s,n$ and $m$; thus, 
there exists $\e>0$ such that
\begin{equation}\label{e28b}
 \sum_{1\le i,j\le \k_{\e}}\left(I^{n,m,i,j}_{1}+I^{n,m,i,j}_{3}\right) \le \sum_{1\le i,j\le \k_{\e}}\frac{\d}{3}
 \PP\left(\mathcal{Q}^{i,j}_{n,\e}\right)\le \frac{\d}{3}
\end{equation}
since $\mathcal{Q}^{i,j}_{n,\e}$, for $1\le i,j\le \k_{\e}$, are pairwise disjoint.

\begin{remark}\label{rr2}
Clearly, $\Sigma_{1}\Sigma_{1}^{\ast}$ is a symmetric $d\times d$ matrix that satisfies the
non-degeneracy condition
\begin{equation}\label{e20bis}
  \langle\Sigma_{1}\Sigma_{1}^{\ast}\x,\x\rangle=|\Sigma_{1}^{\ast}\x|^{2}\ge\l^{2}|\x|^{2},\qquad
  \x\in\R^{d},
\end{equation}
since, by \eqref{e20},
  $$\l|\x|^{2}\le \int_{\R^{N}}\langle\s_{1}(t,x,y)\x,\x\rangle\m(dy)
  =\langle\x,\Sigma^{\ast}_{1}(t,x,\m)\x\rangle
  \le |\x\|\Sigma^{\ast}_{1}(t,x,\m)\x|.$$
\end{remark}
To estimate $I^{n,m,i,j}_{2}$ we use the Krylov's bounds for the process
$(\tilde{X}^{n}_{1,t},\tilde{Y}^{n}_{1,t})$ which, by Lemma \ref{p2} and Remark \ref{rr2}, solves
a non-degenerate SDE. Setting
  $$
  g^{n}_{ij}(t,x_{1},y_{1}):=b^{n}(t,\phi_{i}(t),x_{1},\phi_{j}(t),y_{1}),\qquad
  t\in[0,T],\, x_{1},y_{1}\in\R^{d},$$
we have {\begin{align}
  I^{n,m,i,j}_{2}&\le \int_{0}^{t}\EE\left[\caratt_{\mathcal{Q}^{i,j}_{n,\e}}\left|b^{n}
  (s,\phi_{i}(s),\tilde{X}^{n}_{1,s},\phi_{j}(s),\tilde{Y}^{n}_{1,s})-%
  b^{m}(s,\phi_{i}(s),\tilde{X}^{n}_{1,s},\phi_{j}(s),\tilde{Y}^{n}_{1,s})\right|\right]ds\le
\intertext{(by the Krylov's bounds of Lemma 3 in \cite{MR4421344})}\label{e28c}
  &\le C\|g^{n}_{ij}-g^{m}_{ij}\|_{L^{2d+1}([0,T]\times D_{M_{\e}}\times D_{M_{\e}})}
\end{align}
where
 $D_{M_{\e}}=\{x_{1}\in\R^{d}\mid |x_{1}|<M_{\e}\}$ with $M_\e$ defined in \eqref{e16}
and the positive constant $C$ depends only on $d, M_{\e}$, the linear growth coefficient
$\mathbf{c}$ and the ellipticity constant $\l$}. From \eqref{e28c} it follows that
\begin{equation}\label{e28d}
  \sum_{1\le i,j\le \k_{\e}}I^{n,m,i,j}_{2}\le \frac{\d}{3}\sum_{1\le i,j\le \k_{\e}}\PP\left(\mathcal{Q}^{i,j}_{n,\e}\right)
  \le \frac{\d}{3}
\end{equation}
if $m,n$ are suitably large.

Combining \eqref{e28}, \eqref{e28b} and \eqref{e28d}, for $\e>0$ suitably small and $m\in\N$
suitably large, we have
\begin{equation}\label{e38}
  \lim_{n\to \infty}\mathbf{E}_{1}^{n,m}\le \d.
\end{equation}

\item[$\diamond$] {\bf The term $\mathbf{E}_{3}^{m}$.}\\
The term $\mathbf{E}_{3}^{m}$ can be estimated exactly like $\mathbf{E}_{1}^{n,m}$ after extending
Krylov's bounds to the limiting process
$\tilde{Z}^{\infty}_{1,t}:=(\tilde{X}^{\infty}_{1,t},\tilde{Y}^{\infty}_{1,t})$. Thus, we want to
prove that for every bounded, non-negative and Borel measurable function $g$ with support in
$[0,T]\times D_{M}\times D_{M}$ we have
\begin{equation}\label{e29}
  \int_{0}^{T}\EE\left[g(t,\tilde{Z}^{\infty}_{1,t})\right]dt\le C \|g\|_{L^{2d+1}([0,T]\times D_{M}\times D_{M})}
\end{equation}
with $C$ dependent only on $d, M,\mathbf{c}$ and $\l$. The bound is pretty direct for continuous
functions, indeed if $g$ is also continuous then by weak convergence we have
\begin{equation}\label{e30}
 \int_{0}^{T}\EE\left[g(t,\tilde{Z}^{\infty}_{1,t})\right]dt
 =\lim_{n\to \infty}\int_{0}^{T}\EE\left[g(t,\tilde{Z}^{n}_{1,t})\right]dt\le C \|g\|_{L^{2d+1}([0,T]\times D_{M}\times D_{M})}.
\end{equation}
Next, we prove \eqref{e29} for $g=\caratt_{K}$ where $K$ is a compact subset of $[0,T]\times
D_{M}\times D_{M}$: by Urysohn's lemma, there exists a sequence of continuous functions,  with
values in $[0,1]$, such that $g_{n}\searrow \caratt_{K}$. By Lebesgue's theorem, we have
\begin{align}\label{e31}
 \n(K)&:=\int_{0}^{T}\EE\left[\caratt_{K}(t,\tilde{Z}^{\infty}_{1,t})\right]dt
 =\lim_{n\to \infty}\int_{0}^{T}\EE\left[g_{n}(t,\tilde{Z}^{\infty}_{1,t})\right]dt\le
\intertext{(by \eqref{e30})}
 &\ \le C \lim_{n\to \infty}\|g_{n}\|_{L^{2d+1}([0,T]\times D_{M}\times D_{M})}=C\|g\|_{L^{2d+1}([0,T]\times
D_{M}\times D_{M})}.
\end{align}
Now, $\v$ in \eqref{e31}, as any finite measure on Borel sets, is
regular\footnote{$\v(D)=\sup\{\v(K)\mid K\subseteq D,\ K \text{ compact}\}$ for any Borel set
$D$.}: thus, estimate \eqref{e29} extends to $g=\caratt_{D}$ where $D$ is any Borel subset of
$[0,T]\times D_{M}\times D_{M}$. As a final step, for any non-negative and Borel measurable
function $g$ there exists an increasing sequence $g_{n}$ of simple functions (i.e. linear
combinations of indicator functions of Borel sets) such that $g_{n}(t,z)\nearrow g(t,z)$ for any
$(t,z)$: an easy application of Beppo Levi's theorem allows to conclude the proof of \eqref{e29}.


\item[$\diamond$] {\bf The term $\mathbf{E}_{2}^{n,m}$.}\\
Let $\e$ and $m$ be fixed as previously, ensuring that \eqref{e38} and an analogous estimate for
$\mathbf{E}_{3}^{m}$ hold. Just as for \eqref{e28}, we also have
 $$\mathbf{E}_{2}^{n,m}\le \e\,CT(1+\EE\left[|\y_{x}|^{4}\right])+I^{n,m}$$
where
  $$I^{n,m}:=\int_{0}^{t}\EE\left[\caratt_{\mathcal{Q}_{n,\e}}
  \EE\left[\caratt_{\mathcal{Q}_{n,\e}}\left|b^{m}(s,\x,\tilde{Y}^{n}_{s})-
  b^{m}(s,x,\tilde{Y}^{\infty}_{s})\right|
  \right]\!\!\raisebox{-0.11ex}{$\Big\vert$}_{\x=\tilde{X}^{n}_{s},\,
  x=\tilde{X}^{\infty}_{s}}\right]ds$$
with $\mathcal{Q}_{n,\e}$ as in \eqref{e14Q}. 
Now, $b^{m}(s,\cdot,\cdot)$ is Lipschitz   continuous on $Q_{M_{\e}}$ in \eqref{e25}, uniformly in
$s\in[0,T]$. Thus, we have $I^{n,m}\le I^{n,m}_{1}+I^{n,m}_{2}$ where
\begin{align}
  I^{n,m}_{1}&:= \int_{0}^{t}\EE\left[\caratt_{\mathcal{Q}_{n,\e}}
  \EE\left[\caratt_{\mathcal{Q}_{n,\e}}\left|b^{m}(s,\x,\tilde{Y}^{n}_{s})-
  b^{m}(s,x,\tilde{Y}^{n}_{s})\right|
  \right]\!\!\raisebox{-0.11ex}{$\Big\vert$}_{\x=\tilde{X}^{n}_{s},\,
  x=\tilde{X}^{\infty}_{s}}\right]ds
  \xrightarrow[\ n\to\infty\ ]{}0
\end{align}
by weak convergence of $\tilde{X}^{n}_{s}$ to $\tilde{X}^{\infty}_{s}$ for any fixed $s\in[0,t]$,
and by dominated convergence in the time-integral. Similarly, we have
\begin{align}
  I^{n,m}_{2}&=\int_{0}^{t}\EE\left[\caratt_{\mathcal{Q}_{n,\e}}
  \EE\left[\caratt_{\mathcal{Q}_{n,\e}}\left|b^{m}(s,x,\tilde{Y}^{n}_{s})-
  b^{m}(s,x,\tilde{Y}^{\infty}_{s})\right|
  \right]\!\!\raisebox{-0.11ex}{$\Big\vert$}_{x=\tilde{X}^{\infty}_{s}}\right]ds\xrightarrow[n\to\infty]{}0
\end{align}
by weak convergence of $\tilde{Y}^{n}_{s}$ to $\tilde{Y}^{\infty}_{s}$.
\end{itemize}

\medskip\noindent{\bf Step III: convergence in probability of stochastic integrals.}\\
We show convergence in probability
\begin{equation}\label{e33b}
  \mathbf{I}^{n}:=\int_{0}^{t}\Sigma_{1}^{n}(s,\tilde{X}^{n}_{s},
  \m_{\tilde{Y}^{n}_{s}})d\tilde{W}^{x,n}_{s}-\int_{0}^{t}\Sigma_{1}(s,\tilde{X}^{\infty}_{s},
  \m_{\tilde{Y}^{\infty}_{s}})d\tilde{W}^{x,\infty}_{s}\xrightarrow[n\to\infty]{\PP}0
\end{equation}
by proving that, for any positive constants $c$ and $\d$, we have
\begin{equation}\label{e34}
  \lim\limits_{n\to \infty}\PP\left(\left|\mathbf{I}^{n}\right|>c\right)\le \d.
\end{equation}
As in Step II (cf. \eqref{e33}), for $m\in\N$ we write
$\mathbf{I}^{n}=\mathbf{I}_{1}^{n,m}+\mathbf{I}_{2}^{n,m}+\mathbf{I}_{3}^{m}$ where
\begin{align}
 \mathbf{I}_{1}^{n,m}&:=\int_{0}^{t}\left(\Sigma_{1}^{n}(s,\tilde{X}^{n}_{s},
  \m_{\tilde{Y}^{n}_{s}})-\Sigma_{1}^{m}(s,\tilde{X}^{n}_{s},
  \m_{\tilde{Y}^{n}_{s}})\right)d\tilde{W}^{x,n}_{s},\\
  \mathbf{I}_{2}^{n,m}&:=\int_{0}^{t}\Sigma_{1}^{m}(s,\tilde{X}^{n}_{s},
  \m_{\tilde{Y}^{n}_{s}})d\tilde{W}^{x,n}_{s}-\int_{0}^{t}\Sigma_{1}^{m}(s,\tilde{X}^{\infty}_{s},
  \m_{\tilde{Y}^{\infty}_{s}})d\tilde{W}^{x,\infty}_{s},\\
  \mathbf{I}_{3}^{m}&:=\int_{0}^{t}\left(\Sigma_{1}^{m}(s,\tilde{X}^{\infty}_{s},
  \m_{\tilde{Y}^{\infty}_{s}})-\Sigma_{1}(s,\tilde{X}^{\infty}_{s},
  \m_{\tilde{Y}^{\infty}_{s}})\right)d\tilde{W}^{x,\infty}_{s}.
\end{align}
By Markov inequality and It\^o's isometry, we have
  $$\PP\left(\left|\mathbf{I}_{1}^{n,m}\right|>c\right)\le
  \frac{1}{c^{2}}\EE\left[
  \int_{0}^{t}\left|\Sigma_{1}^{n}(s,\tilde{X}^{n}_{s},
  \m_{\tilde{Y}^{n}_{s}})-\Sigma_{1}^{m}(s,\tilde{X}^{n}_{s},
  \m_{\tilde{Y}^{n}_{s}})\right|^{2}ds\right].$$
Then, proceeding as in the proof of estimate \eqref{e38} for $\mathbf{E}_{1}^{n,m}$ in
\eqref{e36}, we can show that for any $c,\d>0$ there exists a suitably large $m\in\N$ such that
  $$\lim\limits_{n\to \infty}\PP\left(\left|\mathbf{I}_{1}^{n,m}\right|>c\right)\le\d.$$
Similarly, the term $\mathbf{I}_{3}^{m}$ can be treated as $\mathbf{E}_{3}^{m}$ in \eqref{e37}.

Let us assume that $m$ is fixed as stated above. We are left with the term $\mathbf{I}_{2}^{n,m}$,
which sets itself apart from the other two due to the stochastic integrals being associated with
distinct Brownian motions. To deal with this term, we rely on the following lemma, that is a
specific instance of the results in Chapter 2, Section 3 of \cite{MR185620}.
\begin{lemma}[Skorokhod]\label{l3}
On a probability space $(\O,\F,\PP)$, let $W^{n}$ be a sequence of Brownian motions and $f^{n}$ be
a sequence of stochastic processes such that $(f^{n}_{t},W^{n}_{t})$ converges in probability to
$(f^{\infty}_{t},W^{\infty}_{t})$ for any $t\in[0,T]$, where $W^{\infty}$ is a Brownian motion.
Suppose that:
\begin{itemize}
  \item[a)] the stochastic integrals
  $$\int_{0}^{t}f^{n}_{s}dW^{n}_{s},\qquad t\in[0,T],$$
are well-defined for any $n\in\N\cup\{\infty\}$;
  \item[b)] $f^{n}$ are uniformly bounded, i.e. $|f^{n}_{t}(\o)|\le C$ for any $n\in\N$,
  $t\in[0,T]$ and $\o\in\O$;
  \item[c)] for any $c>0$
    $$\lim_{h\to 0}\lim_{n\to\infty}\sup_{0\le t,s\le T\atop |t-s|\le h}\PP(\left|f^{n}_{t}-f^{n}_{s}\right|>c)=0.$$
\end{itemize}
Then
  $$\int_{0}^{t}f^{n}_{s}dW^{n}_{s}\xrightarrow[n\to\infty]{\PP}\int_{0}^{t}f^{\infty}_{s}dW^{\infty}_{s},
  \qquad t\in[0,T].$$
\end{lemma}
We apply Lemma \ref{l3} with
$f^{n}_{t}=\Sigma_{1}^{m}(t,\tilde{X}^{n}_{t},\m_{\tilde{Y}^{n}_{t}})$ and
$W^{n}_{t}=\tilde{W}^{x,n}_{t}$, for $n\in\N\cup\{\infty\}$. Conditions a), b) and c) are readily
verified since $m\in\N$ is fixed and it is clear that $f^{n}$ are uniformly bounded and uniformly
continuous in virtue of the boundedness and smoothness of $\Sigma^{m}_{1}$, and estimate
\eqref{e5Z}. Moreover, we already observed in Step 4 that the stochastic integrals are
well-defined and $\tilde{W}^{x,n}_{t}$ converges in probability to $\tilde{W}^{x,\infty}_{t}$ for
any $t\in[0,T]$, by \eqref{e7}. To conclude, we show convergence in probability (actually, even in
$L^{1}(\O,\PP)$) of $f^{n}_{t}$ to $f^{\infty}_{t}$ for any $t\in[0,T]$. We have
 $$\EE\left[\left|f^{n}_{t}-f^{\infty}_{t}\right|\right]\le
 \mathbf{E}_{1}^{n}+\mathbf{E}_{2}^{n}+\mathbf{E}_{3}^{n}$$
where
\begin{align}
 \mathbf{E}_{1}^{n}&:=\EE\left[\EE\left[\left|\s^{m}_{1}(t,x,\tilde{Y}^{n}_{t})-
  \s_{1}^{m}(t,x,\tilde{Y}^{\infty}_{t})\right|
  \right]\!\!\raisebox{-0.11ex}{$\Big\vert$}_{x=\tilde{X}^{\infty}_{t}}\right]
  \xrightarrow[n\to\infty]{}0
\intertext{by weak convergence of $\tilde{Y}^{n}_{t}$ to $\tilde{Y}^{\infty}_{t}$ in the inner
expectation and Lebesgue dominated convergence theorem in the outer expectation; for any fixed but
arbitrary $\e>0$}
 \mathbf{E}_{2}^{n}&:=\EE\left[\caratt_{\left|\tilde{X}^{n}_{t}-\tilde{X}^{\infty}_{t}\right|\le \e}
  \int_{\R^{N}}\left|\s^{m}_{1}(t,\tilde{X}^{n}_{t},y)-\s^{m}_{1}(t,\tilde{X}^{\infty}_{t},y)\right|
  \m_{\tilde{Y}^{n}_{t}}(dy)\right]\le C_{m}\e
\intertext{where $C_{m}$ is the Lipschitz constant of $\s^{m}_{1}$, and}
  \mathbf{E}_{3}^{n}&:=\EE\left[\caratt_{\left|\tilde{X}^{n}_{t}-\tilde{X}^{\infty}_{t}\right|> \e}
  \int_{\R^{N}}\left|\s^{m}_{1}(t,\tilde{X}^{n}_{t},y)-\s^{m}_{1}(t,\tilde{X}^{\infty}_{t},y)\right|
  \m_{\tilde{Y}^{n}_{t}}(dy)\right]\\
  &\quad\le 2\|\s^{m}_{1}\|_{\infty}\PP\left(\left|\tilde{X}^{n}_{t}-\tilde{X}^{\infty}_{t}\right|> \e\right)
  \xrightarrow[n\to\infty]{}0
\end{align}
by convergence in probability of $\tilde{X}^{n}_{t}$ to $\tilde{X}^{\infty}_{t}$ (cf. \eqref{e7}).
%
%
\endproof

%
%
%

%
%
\subsection{Proof of Lemma \ref{p2}}\label{sp1}
We prove that
\begin{equation}\label{e10}
  \EE\left[\sup_{t\in[0,T]}\left|\tilde{X}^{n}_{t}-\tilde{X}^{n}_{0}-\int_{0}^{t}
  B^{n}(s,\tilde{X}^{n}_{s},\m_{\tilde{Y}^{n}_{s}})ds
  -\int_{0}^{t}\Sigma^{n}(s,\tilde{X}^{n}_{s},\m_{\tilde{Y}^{n}_{s}})d\tilde{W}^{x,n}_{s}\right|\right]=0.
\end{equation}
For simplicity, we consider here only the case $B^{n}\equiv 0$, the general case being completely
analogous. Let $[t]_{m}:=\frac{[2^{m}t]}{2^{m}}$ be the dyadic approximation of $t$; we have
  $$\EE\left[\left|\tilde{X}^{n}_{t}-\tilde{X}^{n}_{0}-\int_{0}^{t}\Sigma^{n}(s,\tilde{X}^{n}_{s},\m_{\tilde{Y}^{n}_{s}})d\tilde{W}^{x,n}_{s}
  \right|^{2}\right]\le E_{1}^{m}+E_{2}^{m}$$
where
\begin{align}
  E_{1}^{m}&:=\EE\left[\left|\tilde{X}^{n}_{t}-\tilde{X}^{n}_{0}-\int_{0}^{t}\Sigma^{n}([s]_{m},\tilde{X}^{n}_{[s]_{m}},\m_{\tilde{Y}^{n}_{[s]_{m}}})d\tilde{W}^{x,n}_{s}
  \right|^{2}\right],\\
  E_{2}^{m}&:=\EE\left[\left|\int_{0}^{t}\left(\Sigma^{n}(s,\tilde{X}^{n}_{s},\m_{\tilde{Y}^{n}_{s}})-
  \Sigma^{n}([s]_{m},\tilde{X}^{n}_{[s]_{m}},\m_{\tilde{Y}^{n}_{[s]_{m}}})\right)d\tilde{W}^{x,n}_{s}
  \right|^{2}\right].
\end{align}
Now, we have
\begin{align}
  E_{1}^{m}&=\EE\left[\left|\tilde{X}^{n}_{t}-\tilde{X}^{n}_{0}-\int_{0}^{t}\Sigma^{n}([s]_{m},\tilde{X}^{n}_{[s]_{m}},\m_{\tilde{Y}^{n}_{[s]_{m}}})d\tilde{W}^{x,n}_{s}
  \right|^{2}\right]\\
  &=\EE\left[\left|\tilde{X}^{n}_{t}-\tilde{X}^{n}_{0}-\sum_{k2^{-m}\le t}
  \Sigma^{n}(k2^{-m},\tilde{X}^{n}_{k2^{-m}},\m_{\tilde{Y}^{n}_{k2^{-m}}})\left(\tilde{W}^{x,n}_{(k+1)2^{-m}}-\tilde{W}^{x,n}_{k2^{-m}}\right)
  \right|^{2}\right]=
\intertext{(by the equivalence of $(\tilde{Z}^{n},\tilde{W}^{x,n})$ and $(Z^{n},W)$)}
  &=\EE\left[\left|X^{n}_{t}-\y-\sum_{k2^{-m}\le t}
  \Sigma^{n}(k2^{-m},X^{n}_{k2^{-m}},\m_{Y^{n}_{k2^{-m}}})\left(W^{n}_{(k+1)2^{-m}}-W^{n}_{k2^{-m}}\right)\right|^{2}\right]\\
  &=\EE\left[\left|{X}^{n}_{t}-\y-\sum_{k2^{-m}\le t}
  \left(\int_{\O}\s^{n}(t,x,{Y}^{n}_{k2^{-m}})d\PP\right)\!\!\raisebox{-0.25ex}{$\Big\vert$}_{x={X}^{n}_{k2^{-m}}}\left({W}_{(k+1)2^{-m}}-{W}_{k2^{-m}}\right)
  \right|^{2}\right]\xrightarrow[\ m\to\infty\ ]{}0,
\end{align}
since $(Z^{n},W)$ solves \eqref{e1b} (with $B^{n}\equiv 0$, by assumption) and by the continuity
and boundedness 
of $\s^{n}$. We also have
  $$\lim_{m\to \infty}E_{2}^{m}=0$$
by standard approximation of a stochastic integral with bounded and continuous coefficient. This
proves that the processes $\tilde{X}^{n}_{t}$ and
  $$\tilde{X}^{n}_{0}+\int_{0}^{t}\Sigma^{n}(s,\tilde{X}^{n}_{s},\m_{\tilde{Y}^{n}_{s}})d\tilde{W}^{x,n}_{s}$$
are modifications, but since they are continuous this suffices to conclude.

\section{Uniqueness: proof of Theorem \ref{at1}}
Let $(X^{i}_{t},W^{i}_{t})_{t\in[0,T]}$, for $i=1,2$, be two weak solutions of the MKV equation
\eqref{e1bb} with flows $\m^{i}=(\m^{i}_{t})_{t\in[0,T]}$ respectively. If we show that
$\m^{1}=\m^{2}$ then weak/strong uniqueness will follow from weak/strong uniqueness of the non-MKV
equation \eqref{e1nonMKV}.

Let us first consider the case where $B_{1}$ is bounded. 
Up to transferring the two solutions to the canonical space (cf., for instance, \cite{MR1780932}
p.152 or \cite{MR1725357}, Theorem IX.1.7.), we can assume that $X^{1},X^{2}$ are defined on the
same filtered probability space and are solutions w.r.t. the same Brownian motion $W$, that is for
$i=1,2$ we have
\begin{equation}\label{eL1}
  \begin{cases}
    dX^{i}_{0,t}=B_{0}(t,X^{i}_{t})dt, \\
    dX^{i}_{1,t}=B_{1}(t,X^{i}_{t},\m_{t}^{i})dt
    +\Sigma_{1}(t,X^{i}_{t})dW_{t}.
  \end{cases}
\end{equation}
Since $B_{0}$ and $\Sigma_{1}$ are independent of $\m$, by Girsanov's theorem we can ``pass'' from one solution
to the other: precisely, let
\begin{align}\label{ge2}
 \l^{\xy}_t &= \Sigma_{1}^{-1}(t,X^1_t)\left(B_1(t,X^1_t,\m^1_t)-B_1(t,X^1_t,\m^2_t)\right),\\ \label{ge3}
 \g^{\xy}_t &= \exp\left({-\int_0^t\l^{\xy}_s dW_s} -\frac{1}{2}\int_0^t|\l^{\xy}_s|^2ds\right),\\
 \label{ge4}
 W^{\xy}_t &= W_t + \int_0^t\l^{\xy}_sds.
\end{align}
Since $B_{1}$ (and consequently $\l^{1}$) is a bounded process, then $\g^{\xy}$ is a martingale and, by Girsanov's theorem, $W^{\xy}$ is a Brownian motion under
the probability measure $P^{\xy}$ defined by $\frac{dP^{\xy}}{dP}\mid_{\F_{t}}=\g^{\xy}_t$.
Therefore, we have
 $$dX^1_{1,t} = B_1(t,X^1_t,\m^2_t)dt + \Sigma_{1}(t,X^1_t)dW^{\xy}_t,$$
and, under the probability $P^{\xy}$, $X^{1}$ is a solution to the same SDE as $X^2$. Due to the
weak uniqueness hypothesis, $X^1$ under $P^{\xy}$ has the same law as $X^2$ under $P$, and in
particular for any Borel set $H\in\mathcal{B}_{N}$ we have
\begin{align}
 \m^{2}_{t}(H) = P(X_{t}^{2}\in H) = P^{\xy}(X_{t}^{1}\in H).
\end{align}
Thus, we can estimate the total variation distance between $\m^{1}_{t}$ and $\m^{2}_{t}$ as
follows:
\begin{align}
 \|\m^1_{t}-\m^2_{t}\|_{\text{TV}}&:=2\sup_{H\in\mathcal{B}_{N}}|\m^1_{t}(H)-\m^2_{t}(H)|\\
 &\ =2\sup_{H\in\mathcal{B}_{N}}|P(X_{t}^{1}\in H)-P^{\xy}(X_{t}^{1}\in H)|\\
 &\ \le\|P\vert_{\F_{t}}-P^{\xy}\vert_{\F_{t}}\|_{\text{TV}}=
\intertext{(by Scheff\'e's theorem (see, for instance, \cite{MR2724359} Sect.2.2.4), since
$\g^{\xy}_t$ is the density of $P^{\xy}$ with respect to $P$ on $\mathcal{F}_t$)}
 &= 2
 \EE\left[1-\g^{\xy}_t\wedge1\right] \le
 2\sqrt{\EE\left[(\g^{\xy}_t)^2\right]-1}\label{gi2}.
\end{align}
Now, by H\"older's inequality we have
\begin{align}
 \EE\left[(\g^{\xy}_T)^2\right]&=
 \EE\left[\exp\left(-2\int_0^T\l^{\xy}_tdW_t-\int_0^T|\l^{\xy}_t|^2dt\right)\right]\\ &\le
 \EE\left[\exp\left(-2\int_0^T\l^{\xy}_tdW_t-4\int_0^T|\l^{\xy}_t|^2dt\right)\exp\left(3\int_0^T|\l^{\xy}_t|^2dt\right)\right]\le
\intertext{(by Cauchy-Bouniakowsky-Schwarz inequality)}
 &\le\sqrt{\EE\left[\exp\left(-4\int_0^T\l^{\xy}_tdW_t-8\int_0^T|\l^{\xy}_t|^2dt\right)\right]}
 \sqrt{\EE\left[\exp\left(6\int_0^T|\l^{\xy}_t|^2dt\right)\right]}\le
\intertext{(since $\exp\left(-4\int_0^T\l^{\xy}_tdW_t-8\int_0^T|\l^{\xy}_t|^2dt\right)$ is
positive local martingale, and therefore a super-martingale)}
 &\le\sqrt{\EE\left[\exp\left(6\int_0^T|\l^{\xy}_t|^2dt\right)\right]}\\
 &=
 \sqrt{\EE\left[\exp\left(6\int_0^T\left|\int_{\R^{N}}
 \Sigma_{1}^{-1}(t,X^1_t)b_1(t,X^1_t,y)\left(\m^1_t(dy)-\m^2_t(dy)\right)\right|^2dt\right)\right]}\\
 \label{gi1}
 &\le
 \sqrt{\EE\left[\exp\left(6\|\Sigma_{1}^{-1}b_1\|_{\infty}^2\int_0^T\|\m^1_t-\m^2_t\|_{\text{TV}}^2dt\right)\right]}.
\end{align}
Notice that, since the last value is non-random, we may drop the expectation: thus, plugging
\eqref{gi1} into \eqref{gi2} and setting
$v(t):=\sup\limits_{s\in[0,t]}\|\m^1_s-\m^2_s\|_{\text{TV}}$, we obtain
 $$v(T)\le 
 2\sqrt{\exp\left(CT v^2(T)\right)-1},
 $$
with {$C:=3\|\Sigma_{1}^{-1}b_1\|_{\infty}^2$}. Let $\a_{0}>0$ be such that $e^{\a}\le 2\a+1$ for
any $0\le \a\le\a_0$: then, since $v(T)\le 2$ by definition,
for $T\le \frac{\a_0}{4C}$ we have
  $$v(T)\le 2\sqrt{2CT}v(T)$$
which implies $v(T)=0$, at least if $2\sqrt{2CT}<1$.

To establish the result for any $T>0$, we proceed via a straightforward inductive argument.
Let us fix some $T>0$ such that $v(T)=0$. We then prove that, for any $k\in\N$, the condition
$v(kT)=0$ implies $v((k+1)T)=0$.
Indeed, let
\begin{align} 
  \gamma_{t}^{(k)}&:= \exp\left(-\int_{kT}^{t}\lambda^1_s dW_s -
  \frac{1}{2}\int_{kT}^{t}|\lambda^1_s|^2ds\right),\\
  W_{t}^{(k)}&:= W_t + \int_{kT}^{t}\lambda^1_s ds,\qquad t\in[kT,(k+1)T].
\end{align}
By Girsanov's theorem $(W_{t}^{(k)})_{t\in[kT,(k+1)T]}$ is a Brownian motion starting at $W_{kT}$
under the probability measure with density $\gamma_{(k+1)T}^{(k)}$ relative to $P$.
Repeating the calculations leading to \eqref{gi2} and \eqref{gi1}, we obtain
  $$v((k+1)T)\leq2\sqrt{\exp\left(CTv^2((k+1)T)\right)-1}$$
with {$C=3\|\Sigma_{1}^{-1}B_1\|_{\infty}^2$}. As earlier, the conditions $2\sqrt{2CT}<1$ and
$T<\frac{\alpha_0}{4C}$ guarantee that $v((k+1)T)=0$ as required. This completes the induction.
\medskip
We now prove the thesis under the 
linear growth condition \eqref{e23c}. We have the following preliminary estimate. Let $X$ be a
weak solution of \eqref{e1bb} on $[0,T]$: there exists a positive constant
$\d=\d(r,\mathbf{c},\l,T)$ such
that
\begin{equation}\label{ge1}
  \EE\left[e^{\d |\bar{X}_{T}|^{2}}\right]<\infty,\qquad \bar{X}_{T}:=\sup_{t\in[0,T]}|X_{t}|.
\end{equation}
Indeed, by \eqref{e23c} we have
\begin{align}
  |X_{t}|&\le|\y|+\mathbf{c}\int_{0}^{t}(1+|X_{s}|])ds+\left|\int_{0}^{t}\Sigma_{1}(s,X_{s})dW_{s}\right|\\
  &\le|\y|+\mathbf{c}T+J_{T}+\mathbf{c}\int_{0}^{t}|X_{s}|ds,
\end{align}
where
  $$J_{T}:=\sup_{t\in[0,T]}\left|\int_{0}^{t}\Sigma_{1}(s,X_{s})dW_{s}\right|,$$
and Gr\"onwall's lemma yields 
\begin{align}\label{ge5}
  \bar{X}_{T}\le
  \left(|\y|+\mathbf{c}T+J_{T}\right)e^{\mathbf{c}T}.
\end{align}
Consequently, we have
\begin{align}\label{ge6}
  \EE\left[e^{\d |\bar{X}_{T}|^{2}}\right]\le \EE\left[e^{3e^{2\mathbf{c}T}\d
  \left(|\y|^{2}+\mathbf{c}^{2}T^{2}+J_{T}^{2}\right)}\right],
\end{align}
which is finite if $\d>0$ is suitably small, thanks to assumption \eqref{e23cd} on the initial
datum and to standard exponential estimates for stochastic integrals, given that $\Sigma_{1}$ is
bounded (cf., instance, Prop. 13.2.4  in \cite{PascucciSC}). This proves \eqref{ge1}. Clearly, if
\eqref{ge1} holds for some $T>0$ then it holds also for any $t\in[0,T]$ with the same $\d$.

Now, let $\l^{1}$ be as in \eqref{ge2}: by \eqref{e23c} we have
\begin{align}
  \EE\left[\exp\left(\frac{1}{2}\int_{0}^{T}|\l^{1}_{t}|^{2}dt\right)\right]&\le
  \EE\left[\exp\left(\mathbf{c}\|\Sigma_{1}^{-1}\|_{\infty}^2\int_{0}^{T}\left(1+|X^{1}_{t}|^{2}\right)dt\right)\right]\\
  &\le
  \EE\left[\exp\left(\mathbf{c}\|\Sigma_{1}^{-1}\|_{\infty}^2
  T\left(1+|\bar{X}^{1}_{T}|^{2}\right)\right)\right]
  <\infty
\end{align}
if $T>0$ is suitably small, thanks to \eqref{ge1}. Thus Novikov's condition is satisfied, $\g^{1}$
in \eqref{ge3} is a martingale and, by Girsanov's theorem, $W^{\xy}$ in \eqref{ge4} is a Brownian
motion under the probability measure $P^{\xy}$ defined by
$\frac{dP^{\xy}}{dP}\mid_{\F_{t}}=\g^{\xy}_t$.

Next, we can retrace the proof given in the case of bounded $b$: using the linear growth condition
\eqref{e23c} and setting $c_{1}=12\mathbf{c}^{2}\|\Sigma_{1}^{-1}\|_{\infty}^2$, estimate
\eqref{gi1} becomes
\begin{align}
 \EE\left[(\g^{\xy}_T)^2\right]&\le
  \sqrt{\EE\left[e^{c_{1}\left(1+|\bar{X}^{1}_{T}|^{2}\right)Tv^{2}(T)}\right]}\\
  \label{ge7}
 &\le e^{\frac{c_{1}}{2}Tv^{2}(T)}\sqrt{\EE\left[e^{c_{1}Tv^{2}(T)|\bar{X}^{1}_{T}|^{2}}\right]}.
\end{align}
Now, recalling that $v(T)\le 2$, as before we have
\begin{equation}\label{ge8}
  e^{\frac{c_{1}}{2}Tv^{2}(T)}\le 1+c_{1}Tv^{2}(T)
\end{equation}
if $T$ is suitably small. On the other hand, we have
\begin{align}
  \EE\left[e^{c_{1}Tv^{2}(T)|\bar{X}^{1}_{T}|^{2}}\right]&=1+\int_{0}^{+\infty}
  2c_{1}Tv^{2}(T)xe^{c_{1}Tv^{2}(T)x^{2}}P(\bar{X}^{1}_{T}\ge x)dx\le
\intertext{(with $\d>0$ as in \eqref{ge1} and using again $v(T)\le 2$)}
 &\le 1+\frac{c_{1}Tv^{2}(T)}{\d}\int_{0}^{+\infty}
  2\d x e^{(4c_{1}T-\d)x^{2}}e^{\d x^{2}}P(\bar{X}^{1}_{T}\ge x)dx\le
\intertext{(if $T$ is sufficiently small so that $4c_{1}T-\d<0$)}
 &\le 1+\frac{c_{1}\EE[e^{\d|\bar{X}^{1}_{T}|^{2}}]}{\d}Tv^{2}(T)
\end{align}
and by the elementary inequality $\sqrt{1+x}\le 1+\frac{x}{2}$ valid for $x\ge-1$,
we have
\begin{equation}\label{ge9}
  \sqrt{\EE\left[e^{c_{1}Tv^{2}(T)|\bar{X}^{1}_{T}|^{2}}\right]}\le
  1+\frac{c_{1}\EE[e^{\d|\bar{X}^{1}_{T}|^{2}}]}{2\d}Tv^{2}(T).
\end{equation}
To sum up, as in \eqref{gi2} we have
\begin{align}
  v(T)&\le 2\sqrt{\EE\left[(\g^{\xy}_t)^2\right]-1}\le
\intertext{(by \eqref{ge7},\eqref{ge8} and \eqref{ge9}, for some positive constant $C$)}
 &\le \sqrt{C T} v(T)
\end{align}
which implies the existence of $T>0$ such that $v(T)=0$. The proof is completed by reiterating the preceding inductive argument.



\begin{footnotesize}
\bibliographystyle{plain}
\bibliography{Bibtex-Final}
\end{footnotesize}

\end{document}